\input amstex
\documentstyle{amsppt}
\input bull-ppt
\NoBlackBoxes
\keyedby{bull253/lbd}
\define\1{\overline}
\define\2{\partial}

\define\vare{\varepsilon}

\define\si{\sigma}

\define\Ga{\Gamma}

\topmatter
\cvol{26}
\cvolyear{1992}
\cmonth{Jan}
\cyear{1992}
\cvolno{1}
\cpgs{119-124}
\title A steepest descent method\\
for oscillatory Riemann-Hilbert problems\endtitle
\author P. Deift and X. Zhou\endauthor
\shorttitle{A steepest descent method}
\address Department of Mathematics, New York 
University--Courant Institute, New
York, New York 10012\endaddress
\address Department of Mathematics, Yale University, New 
Haven, Connecticut
06520\endaddress
\date February 13, 1991 and, in revised form, March 6, 
1991\enddate
\subjclassrev{Primary 35Q20; Secondary 35B40}
\thanks The work of the authors was supported in part by 
NSF Grants DMS-9001857
and DMS-9196033, respectively\endthanks
\endtopmatter

\document
In this announcement we present a general and new 
approach to analyzing the
asymptotics of oscillatory Riemann-Hilbert problems. Such 
problems arise, in
particular, in evaluating the long-time behavior of 
nonlinear wave equations
solvable by the inverse scattering method. We will 
restrict ourselves here
exclusively to the modified Korteweg de Vries (MKdV) 
equation,
$$\gather
y_t-6y^2y_x+y_{xxx}=0,\qquad -\infty<x<\infty,\ t\ge0,\\
\\
y(x,t=0)=y_0(x),\endgather$$
but it will be clear immediately to the reader with some 
experience in the
field, that the method extends naturally and easily to 
the general class of wave
equations solvable by the inverse scattering method, such 
as the KdV, nonlinear
Schr\"odinger (NLS), and Boussinesq equations, etc., and 
also to ``integrable''
ordinary differential equations such as the Painlev\'e 
transcendents.
\par
As described, for example, in \cite{IN} or \cite{BC}, the 
inverse scattering
method for the MKdV equation leads to a Riemann-Hilbert 
factorization problem for
a $2\times 2$ matrix valued function $m=m(\cdot;x,t)$ 
analytic in $\Bbb
C\backslash\bold R$,
$$\left.\matrix m_+(z)=m_-(z)v_{x,t},\qquad z\in\bold R,\\
m(z)\to I\quad\roman{as}\ 
z\to\infty,\endmatrix\right\}\tag1$$
where
$$\gather
m_\pm(z)=\lim_{\vare\downarrow0}m(z\pm i\varepsilon;x,t),\\
\\
v_{x,t}(z)\equiv
e^{-i(4tz^3+xz)\si_3}v(z)e^{i(4tz^3+
xz)\si_3},\qquad\si_3=\left(\matrix
1&0\\0&-1\endmatrix\right),\endgather$$
and
$$v(z)=\left(\matrix 1-|r(z)|^2&-\1{r(z)}\\
r(z)&1\endmatrix\right)=\left(\matrix 
1&-\1r\\0&1\endmatrix\right)
\left(\matrix 1&0\\r&1\endmatrix\right)\equiv b_-^{-1}b_+
.$$
If $y_0(x)$ is in Schwartz space, then so is $r(z)$ and
$$r(z)=-\1{r(-z)},\qquad\sup_{z\in\bold R}|r(z)|<1.$$
From the inverse point of view, given $v(z)$, one
considers a singular integral
\kern-.39pt equation \kern-.39pt (see \kern-.39pt 
\cite{BC}) \kern-.39pt for 
\kern-.39pt the \kern-.39pt associated \kern-.39pt quantity
\kern-.39pt $\mu(z;x\!,t)\!=\!m_+(z;x\!,t)(b_+
^{-1})_{x,t}=m_-(z;x,t)(b_-^{-1})_{x,t}$ and the
solution of the inverse problem is then given by
$$y(x,t)=\left(\left[\si_3,\int_{\bold 
R}\mu(z;x,t)w_{x,t}(z)\frac{dz}{2\pi
i}\right]\right)_{21}\tag2$$
where
$$w_{x,t}=(w_+)_{x,t}+(w_-)_{x,t},\qquad 
w_\pm=\pm(b_\pm-I).$$
\par
Significant work on the long-time behavior of nonlinear 
wave equations solvable
by the inverse scattering method, was first carried out 
by Manakov \cite M and
by Ablowitz and Newell \cite{AN} in 1973. The decisive 
step was taken in 1976
when Zakharov and Manakov \cite{ZM} were able to write 
down precise formulae,
depending explicitly on the initial data, for the leading 
asymptotics for the
KdV, NLS, and sine-Gordon equations, in the physically 
interesting region
$x=O(t)$. A complete description of the leading 
asymptotics of the solution
of the Cauchy problem, with connection formulae between 
different asymptotic
regions, was presented by Ablowitz and Segur \cite{AS}, 
but without precise
information on the phase. The asymptotic formulae of 
Zakharov and Manakov were
rigorously justified and extended to all orders by 
Buslaev and Sukhanov
\cite{BS 1--2} in the case of the KdV equation, and by 
Novokshenov \cite N in
the case of NLS.
\par
The method of Zakharov and Manakov, pursued rigorously in 
\cite{BS} and in
\cite{N}, involves an ansatz for the asymptotic form of 
the solution and
utilizes techniques that are somewhat removed from the 
classical framework of
Riemann-Hilbert problems. In 1981, Its \cite I returned 
to a method first
proposed in 1973 by Manakov in \cite M, which was tied 
more closely to standard
methods for the inverse problem. In \cite I the 
Riemann-Hilbert problem was
conjugated, up to small errors which decay as 
$t\to\infty$, by an appropriate
parametrix, to a simpler Riemann-Hilbert problem, which 
in turn was solved
explicitly by techniques from the theory of isomonodromic 
deformations. This
technique provides a viable, and in principle, rigorous 
approach to the question
of long-time asymptotics for a wide class of nonlinear 
wave equations (see
\cite{IN}). Finally we note that in \cite B, Buslaev 
derived asymptotic
formulae for the KdV equation from an exact determinant 
formula for the
solution of the inverse problem.
\par
In our approach we consider the Riemann-Hilbert problem 
(1) directly, and by
deforming contours in the spirit of the classical method 
of steepest descent, we
show how to extract the leading asymptotics of the MKdV 
equation. In particular
for $x<0$, let $\pm z_0=\pm\sqrt{|x|/12t}$ be the 
stationary phase points for
$i(4tz^3+xz)$. Then the first step in our method is to 
show that (1) can be
deformed to a Riemann-Hilbert problem on a contour 
$\Sigma$ of shape (see
Figure 1), in such a way that the jump matrices $v_{x,t}$ 
on $\bold
R\subset\Sigma$ and on the compact part

\fighere{5pc}
\caption{{\smc Figure 1}}

\noindent of $\Sigma\backslash\bold R$ away from
$\pm z_0$, converge rapidly to the identity as 
$t\to\infty$. Thus we are left
with a Riemann-Hilbert problem on a pair of crosses 
$\Sigma^A\cup\Sigma^B$ (see
Figure 2). As $t\to\infty$, the interaction between 
$\Sigma^A$ and $\Sigma^B$
goes to zero to higher order
and the contribution to $y(x,t)$ through (2) is simply 
the sum of
the contributions from $\Sigma^A$ and $\Sigma^B$ 
separately. Under the scalings
$z\to z(48tz_0)^{-1/2}\mp z_0$, the problems on 
$\Sigma^A$ and $\Sigma^B$ are
reduced to problems on a fixed cross, with jump matrices 
which are independent
of time, and which can be solved explicitly in terms of 
parabolic cylinder
functions, as in \cite I. Substitution in (2) yields, 
finally, the asymptotics
for $y(x,t)$.

\topspace{4pc}
\caption{{\smc Figure 2}}

\par Our result is the following: let
$$\phi(z_0)=\arg\Ga(i\nu)-\frac{\pi}{4}-\arg r(z_0)+
\frac{1}{\pi}\int_{-z_0}
^{z_0}\log|s-z_0|d(\log(1-|r(s)|^2)$$
(here $\Ga$ is the standard gamma function) and let
$$y_a=\left(\frac{\nu}{3tz_0}\right)^{1/2}%
\cos(16tz_0^3-\nu\log(192tz_0^3)+\phi(z_0)),$$
where $\nu=-(2\pi)^{-1}\log(1-|r(z_0)|^2)>0$. Set
$\tau=tz_0^3=(|x|/12t^{1/3})^{3/2}$.
\thm{Theorem} Let $y_0(x)$ lie in Schwartz space with 
reflection coefficient
$r(z)$. As $t\to\infty$, the solution $y(x,t)$ of {\rm 
MKdV} with initial data
$y_0(x)$, has uniform leading asymptotics conveniently 
described at fixed
$t\gg1$, in the six regions shown in Figure $3$.

\fighere{6pc}
\caption{{\smc Figure 3}}

\par
In region {\rm I}, for any $j$,
$$y(x,t)=y_a+O((-x)^{-j}+(-x)^{-3/4}C_j(-x/t))$$
where $C_j(\cdot)$ is rapidly decreasing. In region {\rm 
II},
$$y(x,t)=y_a+(tz_0)^{-1/2}O(\tau^{-1/4}).$$
In region {\rm III},
$$y(x,t)=(3t)^{-1/3}p(x/(3t)^{1/3})+
O(\tau^{2/3}/t^{2/3}),$$
where $p$ is a Painlev\'e function of type {\rm II}.
\par
In region {\rm IV},
$$y(x,t)=(3t)^{-1/3}p(x/(3t)^{1/3})+O(t^{-2/3}).$$
In region {\rm V}, for any $j$,
$$y(x,t)=(3t)^{-1/3}p(x/(3t)^{1/3})+O(t^{-j}+
t^{-2/3}e^{-12\eta\tau^{2/3}})$$
for some $\eta>0$.
\par
Finally, in region {\rm VI}, for any $j$,
$$y(x,t)=O((x+t)^{-j}).\qed$$
\ethm
\rem{Remark $1$} The reader may check that in the overlap 
regions the asymptotic
forms do indeed match. Also the reader may check that in 
regions II and IV, the
formulae for the leading asymptotics agree with those in 
\cite{IN}.
\endrem
\rem{Remark $2$} The above error estimates are not the 
best possible and in region
II in particular, the $\tau^{-1/4}$ decay can certainly 
be improved.
\endrem
\rem{Remark $3$} There is no obstacle in the method to 
obtaining an asymptotic
expansion for $y(x,t)$ to all orders.
\endrem
\rem{Remark $4$} As noted at the beginning of this 
announcement, the method we
have presented extends naturally to the general class of 
nonlinear wave equations
solvable by the inverse scattering
method.
Also, there is no difficulty in
incorporating solutions with solitons.
\endrem

\enddocument